\documentclass[12pt,onecolumn]{IEEEtran}

\usepackage[dvips]{graphicx}
\usepackage[usenames,dvipsnames]{color}
\usepackage{epsfig}
\usepackage{rotating}
\usepackage{amssymb}
\usepackage{amsmath,amsfonts}
\usepackage[boxed]{algorithm}
\usepackage[latin1]{inputenc}
\usepackage{verbatim}
\usepackage[tableposition=top,font=small,labelfont=bf,format=hang]{caption}
\usepackage{booktabs}
\newtheorem{lem}{Lemma}

\newtheorem{thm}{Theorem}
\newtheorem{prop}{Proposition}
\newtheorem{conj}{Conjecture}
\newtheorem{cor}{Corollary}

\title{Extreme values of the Dedekind $\Psi$ function}
\author{ Patrick Sol\'{e}\thanks{Telecom ParisTech,  46 rue Barrault, 75634
Paris Cedex 13, France.}\hspace{1cm}Michel Planat\thanks{ Institut FEMTO-ST, CNRS, 32 Avenue de
l'Observatoire, F-25044 Besan\c con, France}  }

\begin{document}
\maketitle

\begin{abstract}
Let $\Psi(n):=n\prod_{p \vert n}(1+\frac{1}{p})$ denote the Dedekind $\Psi$ function. Define, for $n\ge 3,$ the ratio $R(n):=\frac{\Psi(n)}{n\log\log n}.$
We prove unconditionally that $R(n)< e^\gamma$ for $n\ge 31.$ Let $N_n=2\cdots p_n$ be the primorial of order $n.$ We prove that the statement  $R(N_n)>\frac{e^\gamma}{\zeta(2)}$ for $n\ge 3$ 
is equivalent to the Riemann Hypothesis. 
\end{abstract}

{\bf MsC codes:} 11N56, 11A25, 11M06

\section{Introduction}
The Dedekind $\Psi$ function is an arithmetic multiplicative function defined for every integer $n>0$ by
$$\Psi(n):=n\prod_{p \vert n}(1+\frac{1}{p}).$$
It occurs naturally in questions pertaining to dimension of spaces of modular forms \cite{modular} and to the commutation of operators in quantum physics \cite{P2}.
It is related to the sum of divisor function
$$\sigma(n)=\sum_{d \vert n} d$$
by the inequalities
$$\Psi(n) \le \sigma(n),$$
and the fact that they coincide for $n$ squarefree.
It is also related to Euler $\varphi$ function by the inequalities
$$ n^2>\varphi(n)\Psi(n)> \frac{n^2}{\zeta(2)}  $$
derived in Proposition \ref{trick} below.\\

In view of the studies of large values of $\sigma$ \cite{R} and of low values of $\varphi$ \cite{N}, it is natural to study both the large and low values of 
$\Psi.$ To that end, we define the ratio $R(n):=\frac{\Psi(n)}{n\log\log n}.$ The motivation for this strange quantity is the asymptotics of 
Proposition \ref{mertens}.
We prove unconditionally that $R(n)  < e^\gamma,$ for $n\ge 31$ in Corollary \ref{upper}. Note that this bound would follow also from the
Robin inequality 
$$ \sigma(n)\le e^\gamma n \log \log n $$ for $n\ge 5041$
under Riemann Hypothesis (RH) \cite{R}, since $\Psi(n) \le \sigma(n).$ \\
In the direction of lower bounds, we prove that the statement  $R(N_n)>\frac{e^\gamma}{\zeta(2)}$ for $n\ge 3$ 
is equivalent to RH, where $N_n=2\cdots p_n$ is the primorial of order $n.$ The proof relies on Nicolas's work on the Euler totient function \cite{N}.
%%%%%%%%%%%%%%%%%%%%%%%%%%%%%%%%%%%%%%%%%%%%%%%%%%%%%%%%%%%%
\section{Reduction to primorial numbers}
Define the primorial number $N_n$ of index $n$ as the product of the first $n$ primes
$$N_n=\prod_{k=1}^np_k, $$
so that $N_1=2,\, N_2=6,\cdots$ and so on. 
As in \cite{N}, the primorial numbers play the role here of superabundant numbers in \cite{R}.
They are champion numbers (ie left to right maxima) of the function $x \mapsto \Psi(x)/x:$
\begin{equation}
\frac{\Psi(m)}{m}<\frac{\Psi(n)}{n}~ \mbox{for}~\mbox{any}~  m<n, 
\label{newSuperab} 
\end{equation}
We give  a proof of this fact, which was observed in \cite{P}.
{\prop The primorial numbers $N_n$ are exactly the champion numbers of the function $x \mapsto \Psi(x)/x.$ }
\begin{IEEEproof}
The proof is by induction on $n$. The induction hypothesis $H_n$ is that the statement is true up to $N_n.$ It is clear that $H_2$ is true. Let $N_n \le m <N_{n+1}$ be a generic integer. The number $m$ has at most $n$ distinct prime factors. This, in combination with the observation that $1+1/x$ is monotonically decreasing as a function of $x$, shows 
that $\Psi(m)/m \le \Psi(N_n)/N_n$. Further $\Psi(N_n)/N_n<\Psi(N_{n+1})/N_{n+1}$. The proof of $H_{n+1}$ follows.
% The proof is by induction on $n$. The induction hypothesis $H_n$ is that the statement is true up to $N_n.$
%OEIS sequence $A060735$ begins $2,4,6\dots$ so that $H_2$ is true. Assume $H_n$ true. Let $N_n\le m< N_{n+1}$ denote a generic integer. The prime divisors of $m$ are
%$\le p_n.$ Therefore $\Psi(m)/m \le \Psi(N_n)/N_n$ with equality iff $m$ is a multiple of $N_n.$ Further $\Psi(N_n)/N_n< \Psi(N_{n+1})/N_{n+1}.$ The proof of 
%$H_{n+1}$ follows.
\end{IEEEproof}

In this section we reduce the maximization of $R(n)$ over all integers $n$ to the maximization over primorials.
%First, we deal with multiples of primorials.
{\prop Let $n$ be an integer $\ge 2.$ For any $m$ in the range  $N_n\le m< N_{n+1}$ one has $R(m)\le R(N_n).$}
\begin{IEEEproof}
Like in the preceding proof we have 
$$\Psi(m)/m \le \Psi(N_n)/N_n$$

Since $0<\log \log 6<\log \log N_n \le \log \log m,$ the result follows.
\end{IEEEproof}

%%%%%%%%%%%%%%%%%%%%%%%%%%%%%%%%%%%%%%%%%%%%%%%%%%%%%%%%%%%%%%%%%%%%%%%%%%%%%%%%%%%%%%
\section{$\Psi$ at primorial numbers}
We begin with an easy application of Mertens formula \cite[Th. 429]{HW}.

{\prop \label{mertens} We have, as $n\rightarrow \infty$ $$\lim R(N_n)=\frac{e^\gamma}{\zeta(2)}\approx 1.08.$$}

\begin{IEEEproof}
Writing $1+1/p=(1-1/p^2)/(1-1/p)$ in the definition of $\Psi(n)$
we can combine the Eulerian product for $\zeta(2)$ with Mertens formula
$$\prod_{p\le x}(1-1/p)^{-1}\sim e^\gamma \log(x)$$ to obtain

$$\frac{\Psi(N_n)}{N_n}\sim \frac{e^\gamma}{\zeta(2)} \log(p_n),$$
Now the Prime Number Theorem \cite[Th. 6, Th. 420]{HW} states that $x \sim \theta(x)$ 
for $x$ large. where $\theta(x)$ stands for Chebyshev's first summatory function:
$$\theta(x)= \sum_{p\le x }\log p.$$ This shows
that, taking $x=p_n$ we have
$$p_n\sim \theta(p_n)=\log(N_n).$$
The result follows.
\end{IEEEproof}

This motivates the search for explicit upper bounds on $R(N_n)$ of the form $\frac{e^\gamma}{\zeta(2)}(1+o(1)).$
In that direction we have the following bound.

{\prop \label{fonda} For $n$ large enough to have $p_n\ge 20000,$ that is $n\ge 2263,$ we have $$ \frac{\Psi(N_n)}{N_n} \le \frac{\exp(\gamma+2/p_n)}{\zeta(2)}(\log\log N_n+\frac{1.125}{\log p_n}) $$}

So, armed with this bound, we derive a bound of the form $R(N_n)< e^\gamma$ for $ n\ge A$, with $A$ a constant.\\

{\cor For $n\ge 4,$ we have $R(N_n)< e^\gamma=1.78\cdots$ }

\begin{IEEEproof}
%The proof is computational: a Fincke-Pohst search ( see \cite{FP}) in a sphere of squared radius $7$ about the origin shows the absence of short vectors in that region.
%Use Lemma 6.4 in \cite{CLMS} with $x=p_n$ and $t=2.$

For $p_n\ge 20000,$ we use the preceding proposition. We need to check that 

$$ \exp(2/p_n) (1+\frac{1.125}{\log(p_n) \log\log(N_n)})\le \zeta(2). $$
Since the LHS is a decreasing function of $n$ it is enough to check this inequality for the first $n$ such that $p_n\ge 20000.$
\\ For $5\le p_n\le 20000,$ that is $3\le n\le 2262$ we simply compute $R(N_n),$ and check that it is $< e^\gamma.$

\end{IEEEproof}

We can extend this Corollary to all integers $>30$ by using the reduction of preceding section, combined with some numerical calculations for $30 < n\le N_4.$

{\cor \label{upper} For $n> 30,$ we have $R(n)< e^\gamma.$ }\\

We prepare for the proof of the preceding Proposition by a pair of Lemmas.
First an upper bound on a partial Eulerian product from \cite[(3.30) p.70]{RS}.
{\lem \label{rs} For $x\ge 2,$ we have $$\prod_{p \le x } (1-1/p)^{-1}\le e^{\gamma} (\log x+\frac{1}{\log x})$$}
Next an upper bound on the tail of the Eulerian product for $\zeta(2).$
{\lem \label{us} For $n\ge 2$ we have  $$\prod_{p >p_n } (1-1/p^2)^{-1}\le \exp(2/p_n)$$}

\begin{IEEEproof}
%The proof is computational: a Fincke-Pohst search ( see \cite{FP}) in a sphere of squared radius $7$ about the origin shows the absence of short vectors in that region.
Use Lemma 6.4 in \cite{CLMS} with $x=p_n$ and $t=2.$
\end{IEEEproof}

We are now ready for the proof of Proposition \ref{fonda}.
\begin{IEEEproof}

Write $$\frac{\Psi(N_n)}{N_n} = \prod_{k=1}^n \frac{1-1/{p_k}^2}{1-1/p_k}$$ and use both lemmas to derive

$$ \frac{\Psi(N_n)}{N_n} \le \frac{\exp(\gamma+2/p_n)}{\zeta(2)}(\log p_n+\frac{1}{\log p_n}). $$

Now we get rid of the first $\log$ in the RHS by the bound of \cite[p.206]{R}

$$\log(p_n)<\log \log N_n+ \frac{0.125}{\log p_n}. $$

The result follows.
 
\end{IEEEproof}
%%%%%%%%%%%%%%%%%%%%%%%%%%%%%%%%%%%%%%%%%%%%%%%%%%%%%%%%%%%%%%%%%%%%%%%%%%%%%%%%%%%%%
\section{Lower bounds}

We reduce first to Euler's $\varphi$ function.

{\prop \label{trick} For $n\ge 2$ we have $$n^2 >\varphi(n)\Psi(n)> \frac{n^2}{\zeta(2)} $$}
\begin{IEEEproof}
The first inequality follows at once upon writing
$$\frac{\varphi(n)\Psi(n)}{n^2}=\prod_{p \vert n} (1-1/p^2),$$
a product of finitely many terms $<1.$
Notice for the second inequality that
$$\frac{\varphi(n)\Psi(n)}{n^2}=\prod_{p \vert n} (1-1/p^2)> \prod_{p } (1-1/p^2),$$

an infinite product that is the inverse of the Eulerian product for $\zeta(2).$

\end{IEEEproof}

{\thm \label{lower} Under RH the ratio $R(N_n)$ is $> \frac{e^\gamma}{\zeta(2)}$ for $n\ge 3$. If RH is false, this is still true for infinitely many $n.$ }

\begin{IEEEproof}

Follows by Proposition \ref{trick}, combined with \cite[Theorem 2]{N}.
\end{IEEEproof}

In view of this result and of numerical experiments the natural conjecture is 

{\conj \label{natural} For all  $n\ge 3$ we have  $R(N_n)> \frac{e^\gamma}{\zeta(2)}.$}

The main result of this note is the following.
{\thm \label{main}  Conjecture \ref{natural} is equivalent to RH.}

\begin{IEEEproof}
If RH is true we refer to the first statement of Theorem \ref{lower}. If RH is false we consider the function
$$g(x):=\frac{e^\gamma}{\zeta(2)}\log\theta(x)\prod_{p\le x}(1+1/p)^{-1},$$

 Observing that $\log \theta(p_n)=\log \log N_n,$ we see that
 $g(p_n) <1$ is equivalent to $R(N_n)> \frac{e^\gamma}{\zeta(2)}.$ 
We need to show that there exists an $x_0\ge 3$ such that $g(x_0)>1$ or equivalently $\log g(x_0) >0$. 
% We need to check that for $x$ large enough $g(x)$ can be $>1$ or equivalently 
%$\log g(x) >0.$ 
Using once again the identity $1-1/p^2=(1-1/p)(1+1/p),$ and \cite[Lemma 6.4]{CLMS}, we obtain, upon writing
$$-\log \zeta(2)=\sum_{p\le x}\log(1-1/p^2)+\sum_{p> x}\log(1-1/p^2), $$
the bound

$$\log g(x) \ge \log f(x)-2/x,$$

where $f$ is the function introduced in \cite[Theorem 3]{N}, that is

$$f(x):={e^\gamma}\log\theta(x)\prod_{p\le x}(1-1/p).$$

We know by \cite[Theorem 3 (c)]{N} that, if RH is false, there is a $0<b<1$ such that $\limsup x^{-b}f(x) >0$ and hence $\lim \sup \log f(x)>> \log x $. Since $2/x =o(\log x)$, the result follows.

%Since for $x$ large $2/x <<x^{-b},$ the result follows.
\end{IEEEproof}
%%%%%%%%%%%%%%%%%%%%%%%%%%%%%%%%%%%%%%%%%%%%%%%%%%%%%%%%%%%%%%%%%%%%%%%%%%%%%%%%%%%%%%%%%%%%%%%%
\section{Conclusion}
In this note we have derived upper and lower bounds on the Dedekind $\Psi$ function. We show unconditionally that the function $\Psi(n)$ satisfies 
the Robin inequality. Since $\psi(n)\le \sigma(n)$ this could be proved under $RH$ \cite{R} or by referring to \cite{CLMS}. Of special interest is Conjecture 1 which is shown here to be equivalent to
RH. We hope this new RH criterion will stimulate research on the Dedekind $\Psi$ function.\\

{\bf Acknowledgements:} The authors thank Fabio Anselmi, Pieter Moree, and Jean-Louis Nicolas for helpful discussions.
%%%%%%%%%%%%%%%%%%%%%%%%%%%%%%%%%%%%%%%%%%%%%%%%

\end{document}